\documentclass[a4paper,12pt]{article}
\usepackage[dvipdfmx]{graphicx}
\usepackage{amsthm,amsmath,amssymb}
\usepackage{latexsym}
\usepackage{graphics}
\usepackage{url}
\usepackage{placeins}
\usepackage{hyperref}
\usepackage{color}
\usepackage{bm}
\usepackage{comment}
\usepackage{statex}

\hypersetup{
     colorlinks=true,
     linkcolor=black,
     filecolor=blue,
     citecolor = black,
     urlcolor=cyan,
     }

\definecolor{blue}{rgb}{0,0,0.9}
\definecolor{red}{rgb}{0.9,0,0}
\definecolor{green}{rgb}{0,0.50,0.10}
\definecolor{violet}{rgb}{0.5804,0.0000,0.8275}  

\def\@themcountersep{}

\newtheorem{THEO}{Theorem}[section]
\newtheorem{ALGo}[THEO]{Algorithm}
\newtheorem{CONJ}[THEO]{Conjecture}
\newtheorem{COND}[THEO]{Condition}
\newtheorem{CORO}[THEO]{Corollary}
\newtheorem{DEFI}[THEO]{Definition}
\newtheorem{EXAMP}[THEO]{Example}
\newtheorem{FACT}[THEO]{Fact}
\newtheorem{HYPO}[THEO]{Hypothesis}
\newtheorem{LEMM}[THEO]{Lemma}
\newtheorem{PROB}[THEO]{Problem}
\newtheorem{PROP}[THEO]{Proposition}
\newtheorem{REMA}[THEO]{Remark}
\newcommand{\theo}{\begin{THEO}}
\newcommand{\algo}{\begin{ALGo} \rm}
\newcommand{\cond}{\begin{COND}}
\newcommand{\conj}{\begin{CONJ}}
\newcommand{\coro}{\begin{CORO}}
\newcommand{\defi}{\begin{DEFI} \rm}
\newcommand{\examp}{\begin{EXAMP} \rm}
\newcommand{\fact}{\begin{FACT}}
\newcommand{\hypo}{\begin{HYPO} \rm}
\newcommand{\lemm}{\begin{LEMM}}
\newcommand{\prob}{\begin{PROB} \rm}
\newcommand{\prop}{\begin{PROP}}
\newcommand{\rema}{\begin{REMA} \rm}
\newcommand{\etheo}{\end{THEO}}
\newcommand{\ealgo}{\end{ALGo}}
\newcommand{\econd}{\end{COND}}
\newcommand{\econj}{\end{CONJ}}
\newcommand{\ecoro}{\end{CORO}}
\newcommand{\edefi}{\end{DEFI}}
\newcommand{\eexamp}{\end{EXAMP}}
\newcommand{\efact}{\end{FACT}}
\newcommand{\ehypo}{\end{HYPO}}
\newcommand{\elemm}{\end{LEMM}}
\newcommand{\eprob}{\end{PROB}}
\newcommand{\eprop}{\end{PROP}}
\newcommand{\erema}{\end{REMA}}

\def\0{\mbox{\bf 0}}
\def\1{\mbox{\bf 1}}
\def\2{\mbox{\bf 2}}
\def\3{\mbox{\bf 3}}
\def\4{\mbox{\bf 4}}
\def\5{\mbox{\bf 5}}
\def\6{\mbox{\bf 6}}
\def\7{\mbox{\bf 7}}
\def\8{\mbox{\bf 8}}
\def\9{\mbox{\bf 9}}

\def\s{\mbox{\boldmath $s$}}

\def\x{\mbox{\boldmath $x$}}
\def\y{\mbox{\boldmath $y$}}

\def\A{\mbox{\boldmath $A$}}
\def\B{\mbox{\boldmath $B$}}

\def\P{\mbox{\boldmath $P$}}

\def\GC{\mbox{$\cal G$}}

\def\OC{\mbox{$\cal O$}}

\def\SC{\mbox{$\cal S$}}

\def\Real{\mbox{$\mathbb{R}$}}

\def\s0{\mbox{\scriptsize \boldmath $0$}}

\def\Real{\mathbb{R}}

\begin{document}

\title{ \Large 
The Largest Unsolved QAP Instance Tai256c Can Be Converted into A 256-dimensional 
Simple BQOP with A Single Cardinality Constraint
} 

\author{
\normalsize 
\thanks{NTT DATA Mathematical Systems Inc., Tokyo
 160-00016, Japan 
({\tt fujii@msi.co.jp}).}, \and \normalsize
Sunyoung Kim\thanks{Department of Mathematics, Ewha W. University, Seoul, 
	52 Ewhayeodae-gil, Sudaemoon-gu, Seoul 03760, Korea 
			({\tt skim@ewha.ac.kr}). 
The research was supported  by   NRF 2021-R1A2C1003810.
}, \and \normalsize
Masakazu Kojima\thanks{Department of Industrial and Systems Engineering,
	Chuo University, Tokyo 192-0393, Japan 
  	 ({\tt kojima@is.titech.ac.jp}).
}, \and \normalsize
Hans D. Mittelmann\thanks{School of Mathematical and Statistical Sciences, 
Arizona State University, Tempe, Arizona 85287-1804, U.S.A.
({\tt mittelma@asu.edu}).},  \and \normalsize
Yuji Shinano\thanks{Department of Applied Algorithmic Intelligence Methods (A$^2$IM), 
Zuse Institute Berlin, Takustrasse 7,14195 Berlin, Germany ({\tt shinano@zib.de}).
}
}

\date{\normalsize\today}
\maketitle

\begin{abstract}
\noindent
Tai256c is the largest unsolved quadratic assignment problem  (QAP) instance in 
QAPLIB; a 1.48\% gap remains 
between the best known feasible objective value and lower bound of 
the unknown optimal value. This paper shows that  
the instance can be converted into a 256 dimensional 
binary quadratic optimization problem (BQOP) with a single cardinality constraint which requires 
the sum of the binary variables to be 92.  The converted BQOP is much simpler 
than the original QAP tai256c and it  also inherits some of the symmetry properties. However, it is still 
very difficult to solve.  We present an efficient branch and bound method for 
improving the lower bound effectively.  A new lower bound with 1.36\% gap is also provided.  
\end{abstract}

\bigskip

{\bf Key words. } 
Quadratic assignment problems,  largest unsolved instance, exploiting symmetry, Lagrangian-DNN relaxation,
Newton-bracketing method,
branch-and-bound methods.

\bigskip

{\bf AMS Classification.} 
90C20,  	
90C22,  	
90C25, 	
90C26.  	

\section{Introduction} 

For a positive integer $n$, we let 
$N=\{1,\ldots,n\}$ represent a set of locations and also a set of facilities.
Given $n \times n$ symmetric matrices $\A =[a_{ik}]$ and 
$\B=[b_{j\ell}]$, 
the quadratic assignment problem (QAP) is 
stated as 
\begin{eqnarray}
\zeta^* = \min_{\pi} \sum_{i \in N}^n\sum_{k\in N}^n a_{ik}b_{\pi(i)\pi(k)}, 
\label{eq:QAP0}
\end{eqnarray}
where $a_{ik}$ denotes the flow between facilities $i$ and $k$, $b_{j\ell}$ 
the distance between locations $j$ and $\ell$, and $(\pi(1),\ldots,\pi(n))$ a permutation of 
$1,\ldots,n$ such that 
$\pi(i) = j$ if facility $i$ is assigned to location $j$. 
We assume that 
the distance $b_{jj}$ from $j \in N$ to itself and the flow $a_{ii}$ from $i \in N$ to 
itself are both zero.

The QAP is known as NP-hard in theory, and solving exactly  large scale instances 
({\it e.g.}, $n \geq 40$) is very difficult in practice. 
To obtain an exact optimal solution,  we basically need two types of techniques. 
The first one is for computing approximate optimal solutions and objective values. 
Heuristic methods such as tabu search, genetic method
and simulated annealing have been developed for the QAP 
\cite{CONNOLLY1990,GAMB1997,SKORIN1990,TAILARD1991}. 
Those methods frequently attain near-optimal solutions, which  also happen 
to be the exact optimal solution in some cases. The exactness is, however, not guaranteed in general.
The objective value $\bar{\zeta}$ obtained by those methods serves as 
an upper bound (UB) for the unknown optimal value $\zeta^*$. 
The second technique is 
to provide a lower bound (LB) $\underline{\zeta}$ for $\zeta^*$. If 
$\underline{\zeta}  = \bar{\zeta}$ holds, then we can conclude that 
$\underline{\zeta}  = \zeta^* = \bar{\zeta}$. Various relaxation methods 
\cite{ANSTREICHER2000,GILMORE1962,LAWLER1963,PRendl09,ZHAO1998} have been 
proposed for computing LBs. 
The two techniques mentioned above play as essential tools in the branch 
and bound (BB) method for QAPs 
\cite{ANSTREICHER2002,CLAUSEN1997,GONCLVES2015,LAWLER1963,PARDALOS1997,ROUCAIROL1987}.

In this paper, we focus our attention on the largest unsolved instance tai256c in 
QAPLIB, and show: \vspace{-2mm} 
\begin{itemize}
\item Tai256c can be converted into a $256$-dimensional binary quadratic optimization 
problem (BQOP) with a single cardinality constraint $\sum_{i=1}^{256} x_i = 92$. \vspace{-2mm}  
\item The converted BQOP satisfies a certain nice symmetry property inherited 
from tai256c.  
\vspace{-2mm} 
\item The proposed BB method can improve LBs for the unknown optimal value. 
\vspace{-2mm} 
\end{itemize}
Moreover, we provide a new LB.

In case of the instance tai256c \cite{BURKARD1997}, $n=256$ and both 
matrices $\A$ and $\B$ satisfy certain symmetries. 
Using the symmetry of the matrix $\A$ presented in 
Section 2.1 of \cite{FISCHETTI2012}, we can convert 
QAP~\eqref{eq:QAP0} to the following BQOP with a 
single cardinality constraint:
\begin{eqnarray}
\zeta^* & = & \min \left\{ \x^T\B\x : \x \in \{0,1\}^n \ \mbox{and } \sum_{i=1}^n x_i = 92 \right\}. 
\label{eq:BQOP0}
\end{eqnarray}
We note that the best known feasible solution $\pi^*$ of tai256c 
with the objective value $\bar{\zeta} = 44759294$  
is converted to a feasible solution $\x^*$ of BQOP~\eqref{eq:BQOP0} 
with the same objective value $\bar{\zeta}$ (see \cite{BURKARD1997}).
Also the best known LB,  
$\underline{\zeta} = 44095032$, for tai256c serves as an LB of BQOP~\eqref{eq:BQOP0} 
(see \cite{ANJOSQAPLIB}).  
More details of this conversion is described in Section~\ref{section:conversion}.

Furthermore, the matrix $\B$ satisfies a symmetry such that 
\begin{eqnarray}
	\P^T \B \P = \B \ \mbox{for every } \P \in \GC, \label{eq:symmetry2} 
\end{eqnarray}
where $\GC$ is a group of permutation matrices satisfying 
\begin{eqnarray}
\left. 
\begin{array}{l}
	\P_1 \P_2 \in \GC \ \mbox{if } \P_1 \in \GC \ \mbox{ and }  \P_2  \in \GC, \\
	\P^{-1} = \P^{T} \in \GC  \ \mbox{if } \P \in \GC.
\end{array}
\right\} \label{eq:group0}
\end{eqnarray}
We describe the symmetry including the orbit 
branching~\cite{OSTROWSKI2011,PFETSCH2019} 
in Section~\ref{section:symmetryB}. 

Exploiting symmetries 
of QAPs in their SDP relaxation was discussed in 
\cite{DEKLERK2010a,DEKLERK2012a,PERMENTER2020} (also \cite{BROSCH2022} 
in their DNN relaxation). However, those results are not relevant to 
the subsequent discussion of this paper.

In Section~\ref{section:conversion}, we show how to convert 
QAP ~\eqref{eq:QAP0} into BQOP~\eqref{eq:BQOP0}.
The conversion can be obtained using the result in Section 2.1 
of \cite{FISCHETTI2012}.
The symmetry in the converted BQOP  inherited from QAP ~\eqref{eq:QAP0} is discussed in Section \ref{section:symmetryB}.
In Section~\ref{section:Gurobi}, we present computational results using 
Gurobi Optimizer (version 9.5.2). We show that Gurobi could improve neither of 
the known LB $\underline{\zeta}$ and UB $\bar{\zeta}$. 
In Section~\ref{section:newLowerBound}, we propose a branch-and-bound method 
with the use of the Newton-bracketing (NB) method \cite{KIM2021} 
to prove that the optimal value is not less than a given $\hat{\zeta}$.
Here  $\hat{\zeta}$  is a target LB chosen in the interval of the best known LB $\underline{\zeta}$  
and UB $\bar{\zeta}$ before starting the BB method. 
If we chose $\hat{\zeta}$ to be the best known objective value $\bar{\zeta} = 44759294$, 
then $\bar{\zeta}$ would be proved to be the optimal value. 
In that case, however, the number of nodes which the BB method needs
to generate 
is estimated at more than $4.1$e$14$; hence the target $\hat{\zeta}=\bar{\zeta}$ is too difficult to attain. 
A target LB 
$44150000$ (1.36\% gap) was attained on Mac Studio (20 cpu).

\section{Conversion from QAP~\eqref{eq:QAP0} to BQOP~\eqref{eq:BQOP0}}

\label{section:conversion}

We present a conversion of QAP~\eqref{eq:QAP0} into BQOP~\eqref{eq:BQOP0} using Theorem 1 in Section 2.1 of \cite{FISCHETTI2012}.
Recall that $a_{ii} = 0$ for every $i \in N$. We say that two facilities $i \in N$ 
and $k \in N$ are {\em clones} if 
\begin{eqnarray*}
& & a_{ik} = a_{ki}, \  
a_{ih} = a_{kh} \ \mbox{and } a_{hi} = a_{hk} \ \mbox{for every } h \in N \backslash \{i,k\}. 
\end{eqnarray*}
We write $i \sim k$ if two facilities $i \in N$ and $k \in N$ are clones. Then 
 $\sim$ becomes an equivalence relation. We denote $\theta(i)$ as 
 the class of 
 facilities $k$ which is equivalent to $i$, {\it i.e.},  
 $\theta(i) = \{ k \in N : i \sim k \}$. Define 
\begin{eqnarray*}
M & = & \left\{ \theta(i) : i \in N \right\}, \\ 
\tilde{a}_{\theta(i)\theta(k)} & = & a_{ik} \ \mbox{for every } i\in N  \ \mbox{and } k \in N \ \mbox{with } i \not= k, \ 
\mbox{or equivalently},  \\ 
\tilde{a}_{uv} & = & a_{ik} \ \mbox{if } u=\theta(i), \ v = \theta(k), \ i\in N  \ \mbox{and } k \in N  \ \mbox{with } i \not= k, \\ 
\mu_u & = & \left|\{ i : \theta(i) = u \} \right| \ \mbox{for every } u  \in M.
\end{eqnarray*}
By \cite[Theorem 1]{FISCHETTI2012},  QAP~\eqref{eq:QAP0} is equivalent to 
\begin{eqnarray}
\zeta^* & = & \min \left\{ \sum_{i \in N}\sum_{u \in M} \sum_{j \in N} \sum_{v \in M} 
\tilde{a}_{uv}b_{ij} x_{iu}x_{jv} : 
\begin{array}{l}
x_{iu} \in \{0,1\} \ ((i,u) \in N \times M), \\ 
\sum_{i \in N} x_{iu} = \mu_u \ (u \in M), \\
\sum_{u \in M} x_{iu} = 1 \ (i \in N)
\end{array}
\right\}.  \label{eq:QAP1}   
\end{eqnarray}
We note that $\tilde{a}_{uu}$ may not be zero ($u \in M$). 

In the QAP instance tai256c, the matrix $A$ satisfies that 
\begin{eqnarray*}
& & a_{ik} = 1 \ (1 \leq i, \ k \leq 92, \ i\not=k) \ \mbox{ and } 
a_{ik} = 0 \ (93 \leq i \leq 256 \ \mbox{or } 93 \leq k \leq 256). 
\end{eqnarray*}
It is straightforward to verify that $M$ consists of two classes 
\begin{eqnarray*}
\{1,\ldots,92\} \ \mbox{ and } \{93,\ldots,256\}, 
\end{eqnarray*}
and that 
\begin{eqnarray*}
\widetilde{\A} & = & 
\begin{pmatrix} \tilde{a}_{11} & \tilde{a}_{12} \\ \tilde{a}_{21}& \tilde{a}_{22} \end{pmatrix} 
= \begin{pmatrix} 1 & 0 \\ 0 & 0 \end{pmatrix}.
\end{eqnarray*}
Therefore, QAP~\eqref{eq:QAP1} turns out to be 
\begin{eqnarray}
\zeta^* & = & \min \left\{ \sum_{i \in N}\sum_{u = 1 }^2 \sum_{j \in N} \sum_{v = 1}^2
\tilde{a}_{uv}b_{ij} x_{iu}x_{jv} : 
\begin{array}{l}
x_{iu} \in \{0,1\} \ ((i,u) \in N \times M), \\ 
\sum_{i \in N} x_{iu} = \mu_u \ (u \in M), \\
\sum_{u \in M} x_{iu} = 1 \ (i \in n)
\end{array}
\right\} \nonumber \\
& = & \min \left\{ \sum_{i \in N} \sum_{j \in N} 
b_{ij} x_{i1}x_{j1} : 
\begin{array}{l}
x_{iu} \in \{0,1\} \ ((i,u) \in N\times\{1,2\}), \\ 
\sum_{i \in N} x_{i1} = 92, \\
\sum_{i \in N} x_{i2} = 164, \\
x_{i1} + x_{i2} = 1 \ (i \in N)
\end{array}
\right\} \label{eq:BQP1} \\
& = & \min \left\{ \sum_{i \in N} \sum_{j \in N} 
b_{ij} x_{i1}x_{j1} : 
\begin{array}{l}
x_{i1} \in \{0,1\} \ (i \in N), \\ 
\sum_{i \in N} x_{i1} = 92 
\end{array}
\right\},  \nonumber 
\end{eqnarray}
which coincides with BQOP~\eqref{eq:BQOP0}. 
We note that $x_{i2}$ serves as a slack variable for each binary variable 
$x_{i1}$ $(i \in N)$ in~\eqref{eq:BQP1}.

\section{Symmetry of the matrix $\B$}  

\label{section:symmetryB}

As mentioned in Section 1, the coefficient matrix $\B$  of the objective function of 
BQOP~\eqref{eq:BQOP0} satisfies~\eqref{eq:symmetry2}, where $\GC$ is  
a group of $n \times n$ permutation matrices satisfying~\eqref{eq:group0}. 
Each $n \times n$ permutation matrix $\P$ corresponds uniquely to a permutation 
$\pi$ of $N$ such that 
\begin{eqnarray*}
\P\x = \x_{\pi} 
\ \mbox{for every } \x \in \{0,1\}^n,  
\end{eqnarray*}
where $\x_{\pi} = (x_{\pi(1)},\ldots,x_{\pi(n)})^T$. 
We denote this correspondence $\P \rightarrow \pi$ by $\pi = \varphi(\P)$, which 
forms an isomorphism from the group of all $n \times n$ permutation matrices 
onto the symmetry group $\SC_n$ of the permutations of $N$. 
We identify the group $\GC$ with the permutation group $\{ \varphi(\P) : \P \in \GC \}$ 
since they are isomorphic. 
With this notation, we can rewrite the condition \eqref{eq:symmetry2} as 
\begin{eqnarray}
	\x_{\pi}^T \B \x_{\pi} = \x^T\B\x \ \mbox{for every } \x \in \{0,1\}^n \ \mbox{and } 
	\pi \in \GC. 
	\label{eq:symmetry3}
\end{eqnarray} 

While software called Nauty \cite{MCKAY2010} could be used to compute $\GC$, 
we instead have computed $\GC$ by a simple implicit enumeration of permutation matrices 
$\P$ satisfying~\eqref{eq:symmetry2}, and found: \vspace{-2mm}
\begin{description}
\item{(a) } $\left|\GC\right| = 2048$ \vspace{-2mm}
\item{(b) } The best known feasible solution $\x^*$ of BQOP~\eqref{eq:BQOP0} 
with the objective value $\bar{\zeta}$ is expanded 
to the set of feasible solutions $\{ \P\x^* : \P \in \GC\}$ with the common 
objective value $\bar{\zeta}$, where $\left|\{ \P\x^* : \P \in \GC\}\right| = 1028$; 
$\P_1\x^* = \P_2\x^*$ can occur for distinct $\P_1 \in \GC$ and $\P_2 \in \GC$.\vspace{-2mm}   
\end{description} 

\section{Numerical experiments using Gurobi}

\label{section:Gurobi} 

We report computational results obtained by a general BQOP solver, Gurobi Optimizer (version 9.5.2) \cite{GUROBI}, for tai256c.
We used the formulation \eqref{eq:BQOP0}  with $x_1$ fixed to 1. 
We will see in Section~\ref{section:branching}
that \eqref{eq:BQOP0} has an optimal solution $\x$ with $x_1 = 1$. 
Numerical experiments  were conducted  on  Intel(R) Xeon(R) CPU E7-8880 v4 (2.20GHz) processors using 32 threads with 2TB of RAM.

Computational experiments were performed iteratively with a starting value at each run initialized as
the last updated incumbent solution, 
till the size of the brand-and-bound tree exceeded 2TB memory. 
The Gurobi parameter \verb|Heuristics=0.8| was used in order to accelerate the 
heuristic search at the 
final run. We tested Gurobi parameter \verb|PreQLinearize|, which controls a  linearization of  the model. At the  
root node, \verb|PreQLinearize=1| provided the lower bound 1202900  and 
\verb|PreQLinearize=2|  the lower bound 384430.663. 
Since these values are lower than the lower bound 4117200 obtained with 
\verb|PreQLinearize=-1|, the default value of \verb|PreQLinearize| 
was used.  
Finally,  we obtained 
a lower bound LB = 4181317 and an upper bound UB = 44780594 
with 6 runs. 

\section{A branch and bound method for a given target lower bound}

\label{section:newLowerBound}

The optimal value $\zeta^*$ of BQOP~\eqref{eq:BQOP0} is currently unknown.
Only a UB $\bar{\zeta} = 44759294 \geq \zeta^*$ and an LB
$\underline{\zeta} = 44095032 \ ($with $1.48\% \ \mbox{gap from } \bar{\zeta}) \ \leq \zeta^*$ are known.
To compute the optimal value $\zeta^*$, we need to improve  
the UB and/or the LB.

Focusing on improving the LB, 
we present a BB  method with the use of the Newton-bracketing (NB) 
method \cite{KIM2021} for solving a Lagrangian doubly nonnegative (Lag-DNN) 
relaxation of subproblems of   
BQOP~\eqref{eq:BQOP0}. The Lag-DNN relaxation is known to provide 
LBs with higher quality for BQOPs than the LP and the SDP relaxation \cite{ITO2017,KIM2013}. 
Before starting the BB method, a target LB, $\hat{\zeta}$, is first set such that 
$\underline{\zeta} = 44095032 < \hat{\zeta} \leq \bar{\zeta} = 44759294$. 
A target LB, $\hat{\zeta}$, is the desired value to attain. 
Ideally, we want set $\hat{\zeta} = \bar{\zeta}$ to confirm whether  $\bar{\zeta}$ 
is the optimal value. But such a setting may be too ambitious, 
which requires much stronger computing power than
the machine currently used.  
As a larger $\hat{\zeta}$ is set,
the computational cost rapidly increases as we will see in 
Figure~1. 

In Section~\ref{section:subproblms}, we describe a class of 
subproblems of BQOP~\eqref{eq:BQOP0} which appear in the enumeration tree generated by the BB method.
For the lower bounding procedure, we use 
the Lag-DNN relaxation of a subproblem 
and the NB method for computing 
its optimal value which serves as an LB of the subproblem. 
We refer to \cite{KIM2021} for details on them. Any upper bounding procedure is not 
incorporated. 
The branching procedure used in the BB method is presented in 
Section~\ref{section:branching}, and 
numerical results in Section~\ref{section:result}. 

\rema The BB method with a fixed target LB above originally developed 
for large scale QAPs was successful to obtain improved lower 
bounds for some of the QAP instances in QAPLIB 
including sko100a,$\ldots, $ sko100f, tai80b, tai100b and tai150b. 
See \cite{ANJOSQAPLIB}. 
The size of QAP tai256c, however, was too large to handle by the original BB method 
for the QAPs. 
\erema

\subsection{A class of subproblems of BQOP~\eqref{eq:BQOP0}} 

\label{section:subproblms}

Let 
\begin{eqnarray*}
\SC & = & \left\{ (I_0,I_1,F) : 
\begin{array}{l}
\mbox{a partition of $N$, {\it i.e.},} \
I_0 \bigcup I_1\bigcup F = N, \\ 
I_0, \ I_1 \ \mbox{and } F \ \mbox{are disjoint with each other}
\end{array}
\right\}. 
\end{eqnarray*}
Obviously, $F$ is uniquely determined as 
$F = N \backslash \big(I_0 \bigcup I_1 \big)$ for each $(I_0,I_1,F) \in \SC$. Hence, 
$F$ in the triplet $(I_0,I_1,F)$ is redundant, and we frequently omit $F$ for the simplicity 
of notation. 
For each $(I_0,I_1,F) \in \SC$, 
we consider a subproblem of BQOP~\eqref{eq:BQOP0} 
\begin{eqnarray*}
\mbox{BQOP}(I_0,I_1) \ : \ \zeta(I_0,I_1) & = & \min
\left\{ \x^T \B \x : 
\begin{array}{l}
\x \in \{0,1\}^n, \ \displaystyle \sum_{i=1}^n x_i = 92, \\
x_i = 0 \ (i \in I_0), \ x_j = 1 \ (j \in I_1)
\end{array} 
\right\}\\
& = & \min \left\{ \y^T \B(I_0,I_1) \y : 
\begin{array}{l}
\y \in \{0,1\}^{F}, \\
 \displaystyle \sum_{i \in F} y_i = 92 - \left| I_1 \right|
\end{array} 
\right\}, 
\end{eqnarray*}
where 
\begin{eqnarray*}
& & \y \in \Real^F \ \mbox{denotes the subvector of $\x$ with elements $x_i$ $(i\in F)$}, \\ 
& & \B(I_0,I_1) = \B_{FF} + 
2 \times \mbox{diagonal matrix of} \left(\sum_{k\in I_1} \B_{kF} \right), \\ 
& & \B_{EF} = \mbox{the $\left|E\right| \times \left|F\right|$ submatrix of $\B$} \
\mbox{consisting of elements $B_{ij}$ $(i \in E, \ j \in F)$}.
\end{eqnarray*}

Applying a penalty function method, we can convert $\mbox{BQOP}(I_0,I_1)$ into a simple 
quadratic unconstrained binary optimization (QUBO) 
\begin{eqnarray*}
\mbox{QUBO}(I_0,I_1,\lambda) \ : \ \zeta(I_0,I_1,\lambda) 
& = & \min \left\{ \y^T \B(I_0,I_1,\lambda) \y : 
\y \in \{0,1\}^{F} 
\right\}, 
\end{eqnarray*}
where 
\begin{eqnarray*}
\y^T \B(I_0,I_1,\lambda) \y 
& = & \y^T \B(I_0,I_1) \y + \lambda \big(\sum_{i \in F} y_i - 92 + \left| I_1 \right|\big)^2 \
\mbox{for every } \y \in \{0,1\}^F, 
\end{eqnarray*}
and $\lambda > 0$ is a Lagrangian or penalty parameter. 
It is straightforward to show that 
$\zeta(I_0,I_1,\lambda)$ converges monotonically to $\zeta(I_0,I_1)$ from below as 
$\lambda \rightarrow \infty$. For computing an LB of $\mbox{BQOP}(I_0,I_1)$ in the BB method 
presented in this section,  
we applied the NB method to the DNN relaxation of $\mbox{BQOP}(I_0,I_1,\lambda)$ with 
$\lambda = 1.0$e$8 / \parallel \B(I_0,I_1) \parallel$. 
See \cite{KIM2021} for more details. 
We note that the matrix $\B(I_0,I_1,\lambda)$ satisfies the 
same symmetry property~\eqref{eq:symmetry2} presented in 
Section~\ref{section:branching}. 
In particular, BQOP~\eqref{eq:BQOP0} 
is converted into QUBO$(\emptyset,\emptyset,\lambda)$.

\subsection{The orbit branching} 

\label{section:branching}

We discuss the orbit branching in ~\cite{OSTROWSKI2011,PFETSCH2019}.
As  mentioned in Section \ref{section:symmetryB}, $\B = \B(\emptyset,\emptyset)$ satisfies the 
symmetry property~\eqref{eq:symmetry3}. This property is partially inherited to 
many $\B(I_0,I_1)$ $((I_0,I_1,F) \in \SC)$. Let $(I_0,I_1,F) \in \SC$ be fixed.
Assume in general that 
\begin{eqnarray}
\y_{\pi}^T \B(I_0,I_1) \y_{\pi} = \y^T \B(I_0,I_1) \y \ 
\mbox{for every } \y \in \{0,1\}^{F} \ \mbox{and } \ \pi \in \GC(I_0,I_1)  
\label{eq:symmetry52}
\end{eqnarray}
holds, where $\GC(I_0,I_1)$ is a permutation group. 
Then, we can define the equivalence relation $\sim$ on $F$ by 
$i \sim j$ if $j = \pi_i$ for some $\pi \in \GC(I_0,I_1)$. Let 
$\omega(i) = \{j \in F: j \sim i \}$ for every $i \in F$, and $\OC(I_0,I_1) = 
\left\{\omega(i): i \in F \right\}$. 
Each $o \in \OC(I_0,I_1)$ is called an {\em orbit} of the permutation 
group $\GC(I_0,I_1)$. We note that $\omega(i) = \omega(j)$ if $i ~ j$, and that
$\OC(I_0,I_1)$ can be the trivial identity permutation $\pi$ such that $\y_{\pi} = \y$ 
for every $\y \in \{0,1\}^{F}$; hence $\OC(I_0,I_1)$ can be defined consistently 
even if $\B(I_0,I_1)$ does not satisfy any symmetry. 
Let $\min(o)$ denote the minimum index of orbit~$o$, which 
serves as a representative from~$o$. 

Assume that $o \in \OC(I_0,I_1)$. 
Obviously, the feasible region of BQOP$(I_0,I_1)$ is 
partitioned disjointly into the subset 
\begin{eqnarray*}
	\left\{ \y \in \{0,1\}^{F} 
	: \displaystyle \sum_{i \in F} y_i = 92 - \left| I_1 \right|, 
	\ y_j = 0 \ \mbox{for every } j \in o \right\},  
\end{eqnarray*}
which forms the feasible region of BQOP$(I_0\bigcup o,I_1)$, and the subset 
\begin{eqnarray*}
	\left\{ \y \in \{0,1\}^{F} 
	: \displaystyle \sum_{i \in F} y_i = 92 - \left| I_1 \right|, \ y_j = 1 
	\ \mbox{for some } j \in o \right\},   
\end{eqnarray*}
which forms the feasible region of BQOP$(I_0,I_1 \bigcup \{j\})$. 
We also see that 
\begin{eqnarray*}
\min \{ \y^T\B(I_0,I_1)\y : \y \in \ \mbox{the latter subset} \}
\end{eqnarray*}
is equivalent to $\min_{j \in o} \zeta(I_0,I_1\bigcup \{j\})$.  
By the construction of orbit $o$, we know that all BQOP$(I_0,I_1\bigcup \{j\})$ 
$(j \in o)$ are equivalent in the sense that they share a common optimal value 
$\zeta(I_0,I_1\bigcup \min (o))$. 
Therefore, we can branch BQOP$(I_0,I_1)$ into two sub BQOPs, 
BQOP$(I_0\bigcup o,I_1)$ and BQOP$(I_0,I_1\bigcup \min(o))$. 

In general, $\OC(I_0,I_1)$ consists of multiple orbits. How  $o$ 
is chosen from $\OC(I_0,I_1)$ for branching of BQOP$(I_0,I_1)$ into 
BQOP$(I_0\bigcup o,I_1)$ and BQOP$(I_0,I_1\bigcup \min(o))$
is an important issue to design an efficient branch and bound method. In our 
numerical experiment presented in Section~\ref{section:result}, \vspace{-2mm} 
\begin{itemize}
\item an orbit $o$ is chosen from $\OC(I_0,I_1)$ 
according to the average objective value of 
BQOP$(I_0,$ $I_1\bigcup\min(o))$ over all 
feasible solutions, so that the chosen orbit, $o^*$,  attains the largest value.  
Then, we apply the branching of BQOP$(I_0,I_1)$ into 
two subproblems BQOP$(I_0$ $\bigcup o^*,I_1)$ and BQOP$(I_0,I_1\bigcup \{\min(o^*)\})$. 
Here the average objective value of BQOP$(I_0,I_1)$ over all 
feasible solutions is computed as the objective value $\x^T\B\x$ with 
$x_i = 0 \ (i\in I_0), \ x_j = 1 \ (j \in I_1)$ and $x_k = (92-|I_1|) / |F| \ (k \in F)$ 
for every $(I_0,I_1,F) \in \SC$.  
\end{itemize}
See \cite[Section 5]{FUJII21} for some other branching rules which can be combined with the NB method. 

$\GC(\emptyset,\emptyset) = \GC$  has the single orbit 
$o = N = \{1,\ldots,n\}$. We branch BQOP$(\emptyset,\emptyset)$ 
into two subproblems BQOP$(N,\emptyset)$ and BQOP$(\emptyset,\{1\})$. 
Obviously, the former BQOP$(N,\emptyset)$ is infeasible. 
Table 1 
summarizes the branching of the node BQOP$(\emptyset,\{1\})$ into 
BQOP$(\{2,16,17,$ $241\},\{1\})$ and BQOP$(\emptyset,\{1,2\})$, where orbit $\{2,16,17,241\}$ 
is chosen from $\OC(\emptyset,\{1\})$. 

\begin{table}[htp]
\scriptsize{
\begin{center}
\caption{A summary of branching of BQOP$(I_0,I_1)$ 
with $I_0 = \emptyset$, $I_1 = \{1\}$ and $F = \{2,3,\ldots,256\}$ 
into BQOP$(\{2,16,17,241\},\{1\})$ and BQOP$(\emptyset,\{1,2\})$. 
Here $F = \{2,3,\ldots,256\}$ is partitioned into 44 orbits, 
which consist of $21$ orbits with size $8$, 
$21$ orbits with size $4$, $1$ orbit with size $2$ and $1$ orbit with size $1$.  
The $44$ orbits are listed according to the decreasing order of the average 
objective value of BQOP$(\emptyset,\{1,\min(o)\})$ over all feasible solutions.
\vspace{3mm} } 
\label{table:orbits}
\begin{tabular}{|c|l|c|c|c|c|}
\hline
Orbit  &       & the size & the average objective \\
number & \multicolumn{1}{c|}{orbit} & of orbit & value of BQOP$(\emptyset,\{1,\min(o)\})$ \\
\hline
 1   & 2 16 17 241    &  4  & 52655297.0 \\ 
 2   & 18 32 242 256    &  4  & 52567852.0 \\
 3   & 3 15 33 225    &  4  & 52524130.0 \\
 4   & 19 31 34 48 226 240 243 255    &  8  & 52515385.0 \\
 5   & 35 47 227 239    &  4  & 52502268.0 \\
     & $\cdots$ & & $\cdots$ \\ 
30   & 87 91 102 108 166 172 183 187    &  8  & 52483274.0 \\
31   & 9 129    &  2  & 52483139.0 \\
32   & 72 74 117 125 149 157 200 202    &  8  & 52483097.0 \\
33   & 25 130 144 249    &  4  & 52483097.0 \\
    & $\cdots$ & & $\cdots$ \\  
43   & 121 136 138 153    &  4  & 52481955.0 \\
44   & 137    &  1  & 52481773.0 \\
\hline
\end{tabular}
\end{center}
}
\end{table}

In addition to the branching rule mentioned above, we employ 
the simple breadth first search; the method to search the enumeration tree is not relevant 
to the computational efficiency since the incumbent objective value is fixed to 
the target LB $\hat{\zeta}$ and 
any upper bounding procedure is not applied.
At each node BQOP$(I_0,I_1)$ of the enumeration tree, 
the NB method generates a sequence of intervals $[a_p,b_p]$ $(p=1,2,\ldots)$ 
satisfying the monotonicity: 
(1) $a_p$  converges monotonically to an LB $\nu$ of BQOP$(I_0,I_1)$
from below, and (2) $b_p$  converges monotonically to $\nu$ from above. 
Thus, if $\hat{\zeta} \leq a_q$ holds for some $q$, we know that the LB $\nu$ to which 
the interval $[a_p,b_p]$ converges is not smaller than $\hat{\zeta}$, 
and BQOP$(I_0,I_1)$ can be pruned. 
On the other hand, if $b_q < \hat{\zeta}$ holds for some $q$, we know the LB $\nu$ 
of BQOP$(I_0,I_1)$ is 
smaller than $\hat{\zeta}$; hence  the iteration can be stopped and  branching  
to BQOP$(I_0,I_1)$ can be applied. Therefore, the above properties (1) and (2) 
of the NB method 
work very effectively to increase the computational efficiency of the BB method. See Figure 3.

\subsection{Numerical results}

\label{section:result}

All the computations for numerical results reported in this section were performed using MATLAB
2022a on Mac Studio with Apple M1 Ultra CPU, 20 cores and 128 GB memory.

To choose a reasonable target LB  \ $\hat{\zeta}$ which can be  attained, we performed 
some preliminary numerical experiments to estimate the computational work. 
Given a target LB $\hat{\zeta}$, 
we construct the enumeration tree by the breadth first search as long as the number $t_k$ 
of nodes at the depth k of the tree is smaller than $1000$. 
Suppose that $t_0, t_1,\ldots,t_{\ell-1} < 1000 \leq t_{\ell}$; 
hence the full enumeration tree has been constructed up to the depth $\ell$ 
by the BB method. We start sampling at the depth $\ell$ and construct a random subtree to estimate the total number 
of nodes of the full enumeration tree.
Let $\bar{t}_{\ell} = t_{\ell}$. 
At each depth $k \geq \ell$, we choose $s_k$ nodes randomly from $\bar{t}_k$ active nodes 
for the next depth $(k+1)$, where 
\begin{eqnarray*}
	s_k & = & \left\{
	\begin{array}{ll}
	100 & \mbox{if } \bar{t}_{k} \geq 500, \\
	\bar{t}_{k} & \mbox{otherwise}.  
	\end{array}
\right. 
\end{eqnarray*}
Then, we apply the lower bounding procedure using the NB method to the selected $s_k$ 
nodes and the branching procedure to the resulting $r_k$ active nodes to generate a subset 
of the nodes of the full enumeration tree at the depth $(k+1)$. 
Next, we let $\bar{t}_{k+1} = 2 * r_k$, 
which is the cardinality of the subset (the number of nodes in the subset) 
as each active node is branched into two child nodes. 
We continue this process till $r_k$ attains $0$. We may regard $2 *r_k /s_k = \bar{t}_{k+1}/s_k$ 
as the increasing rate of the nodes from the depth $k$ to the depth ${k+1}$, and the total 
number of nodes of the full enumeration tree is estimated by 
\begin{eqnarray}
& & \sum_{k=1}^{\ell} t_k  + \sum_{k > \ell} \hat{t}_{k}, \ \mbox{where } \
    \hat{t}_{\ell} = t_{\ell}, \ \hat{t}_{k+1} = (2r_k/s_k)\hat{t}_k \ (k \geq \ell). 
\label{eq:noOfNodes}
\end{eqnarray}
Although this unrefined method seems simple, 
it provides useful information on whether a given target LB 
can be attained by the BB method on the computer used. 

Figure~1 (A) illustrates how the estimated increasing rate changes in 
the eight  cases where the target LB $\hat{\zeta} = 44759294$ (the best known objective 
value $\bar{\zeta}$), 
$44600000$ (with 0.36\% gap from $\bar{\zeta}$), $44400000$ (0.80\% gap), $44200000$ (1.25\% gap), 
$44150000$ (1.36\% gap),  $44130000$ (1.41\% gap), $44120000$ (1.43\% gap) 
and $44100000$ (1.47\% gap) are taken. By applying the formula \eqref{eq:noOfNodes}, 
the total number of nodes in the 
estimated tree is computed as in Table~\ref{table:noOfNodes}. 
The mean execution time to process one node, 
which mainly consists of execution time to solve the Lag-DNN relaxation of 
B$(I_0,I_1)$ by the NB method, is about $50 \sim 60$ seconds, depending on the target LB.  
Obviously, the cases with 
$\hat{\zeta} = 44759294$ and $44600000$ are very challenging. 
The cases with $\hat{\zeta} =44200000$ and $44150000$ 
could be easily  processed by a supercomputer or a large scale clustered computer system. 
We show more details of the numerical results for 
$\hat{\zeta} = 44150000, \ 44130000, \ 44120000,\  44100000 $ cases below.  

Figure 1 (B) compares the increasing rate of nodes of the full enumeration 
tree with its rough estimation described above  for the four cases. Table~\ref{table:noOfNodes} shows 
the total number of nodes of the enumeration tree in these cases.  
Figure 2 (A) displays how the number of nodes 
changes as the depth $k$ increases, and 
Figure 2 (B) how the number of nodes with size 2 orbit changes. All other nodes are of 
the trivial single orbit $N$, except the root node having size 256 orbit as shown in Section 5.2 
and the depth 1 node having size 4 orbit as  observed in Table~\ref{table:orbits}.  

Recall that the NB method applied to a node BQOP$(I_0,I_1)$ generates a sequence of 
interval $[a_p,b_p]$ $(p=1,\ldots)$ which monotonically converges to an LB, $\nu$, of 
BQOP$(I_0,I_1)$. Hence, the iteration terminates either when $b_q < \hat{\zeta}$ occurs 
--- BQOP$(I_0,I_1)$ is turned out to be active in this case --- or when $\hat{\zeta} \leq a_q$ 
occurs --- BQOP$(I_0,I_1)$ is pruned in this case. We see from Figure 3 that 
any tight LB is not necessary to decide whether BQOP$(I_0,I_1)$ is active or to be pruned 
in most cases, particularly, in earlier stage of the BB method. This is an important 
feature of the NB method, which works effectively to increase the computational efficiency,  
when it is incorporated in the BB method.

\begin{figure}
\begin{center}
\caption{(A) Estimated increasing rate of nodes of the full enumeration tree. 
(B) Comparison of increasing rate of nodes of the full enumeration tree and its estimation.} 
\includegraphics[width=0.45\textwidth]{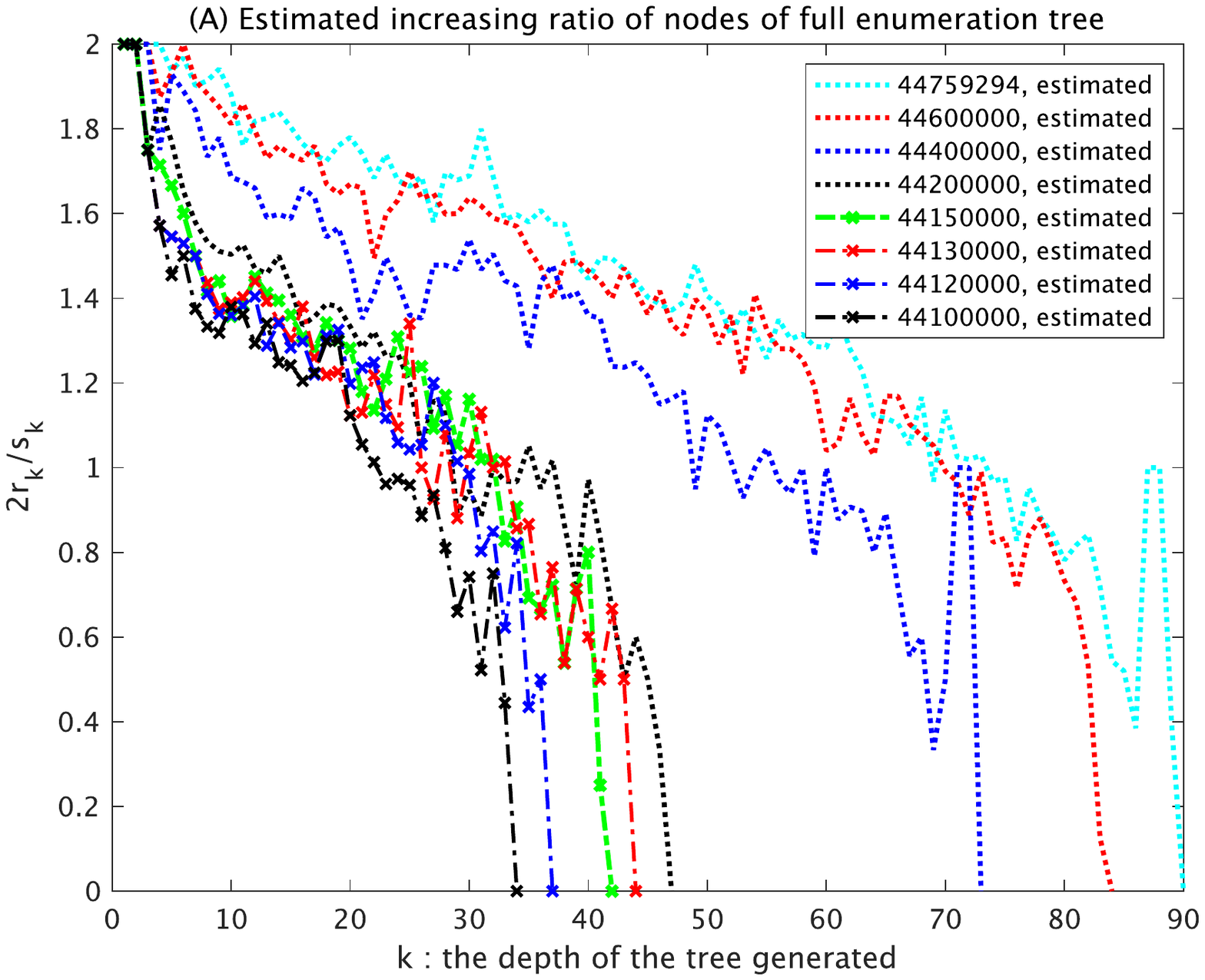}
\includegraphics[width=0.45\textwidth]{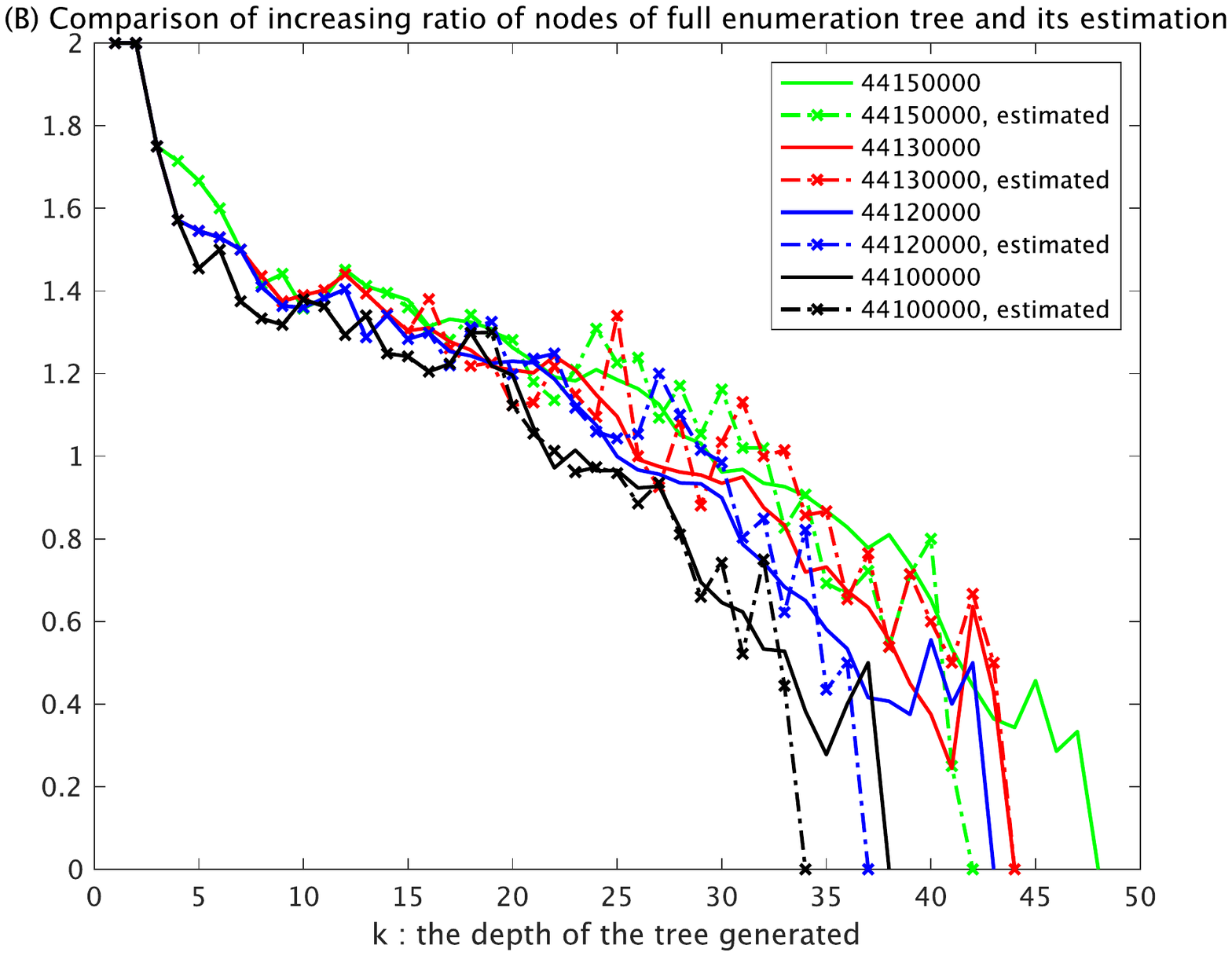}
\end{center}
\end{figure}

\begin{table}[htp]
\scriptsize{
\begin{center}
\caption{Total number (B) of nodes  of the full enumeration tree and its estimation (A).}
\vspace{3mm}
\label{table:noOfNodes}
\begin{tabular}{|c||l|l|l|l|l|l|l|l|}
\hline
 & \multicolumn{8}{|c|}{Target LB} \\ 
                & 44759294 & 44600000 &  44400000 & 44200000 & 44150000 &44130000 & 44120000 & 44100000 \\ 
\hline
(A) Estimated tree  & 4.1e14  & 2.8e13  & 1.4e10   & 9.2e5   & 3.2e5   & 1.2e5  & 9.2e4   & 2.2e4   \\
\hline
(B) Full tree       &          &          &           &          & 277304  & 102310  & 63554    & 23510 \\
\hline 
\end{tabular}
\end{center}
}
\end{table}

\begin{figure}
\begin{center}
\caption{
(A) 
The number of nodes of the enumeration tree at the depth $k$. 
(B) 
The number of nodes of the enumeration tree with size 2 orbit at the depth $k$. 
}
\vspace{-25mm}
\includegraphics[width=0.45\textwidth]{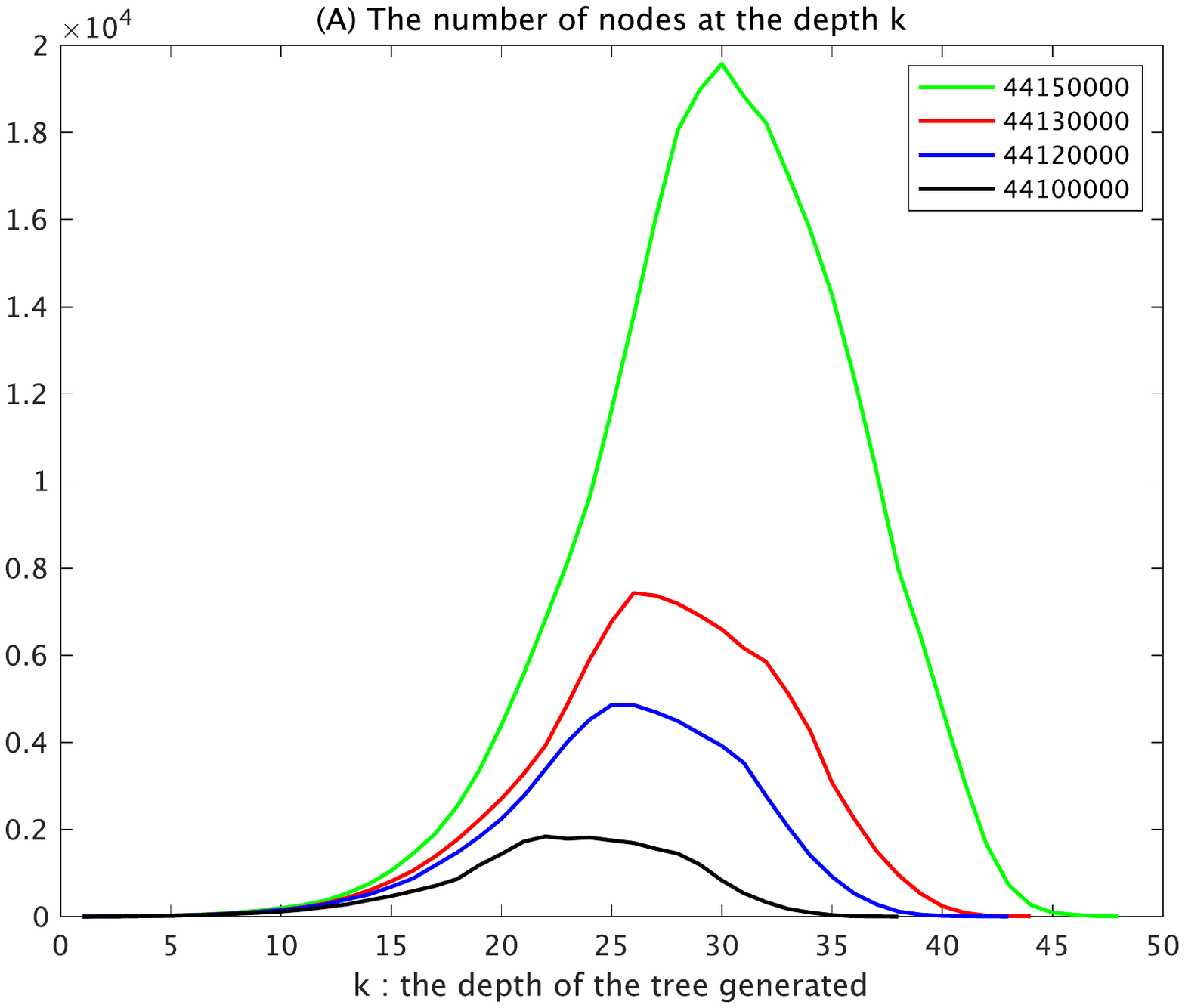}
\includegraphics[width=0.45\textwidth]{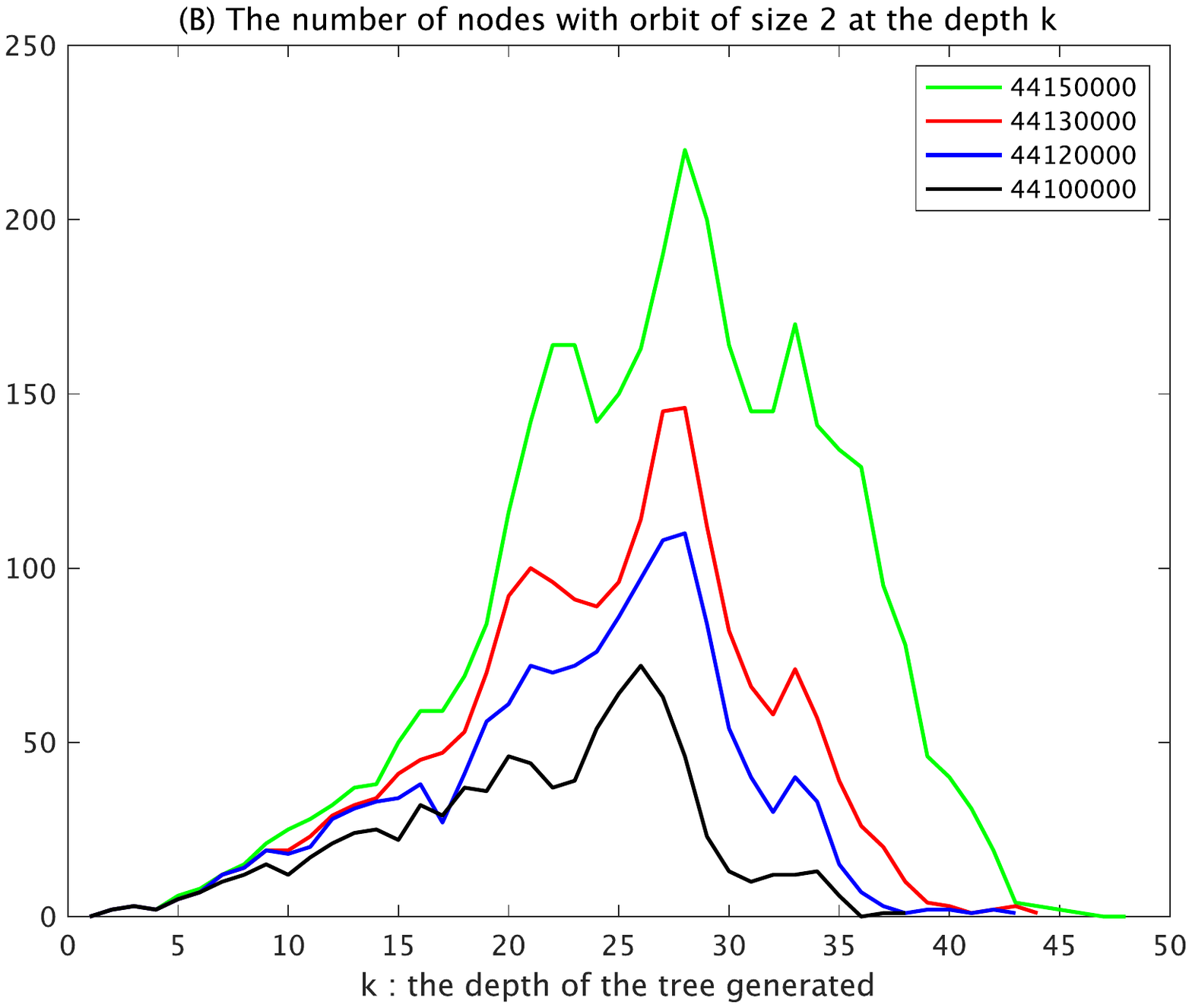}
\end{center}
\end{figure}


\begin{figure}
\begin{center}
\caption{The mean of $a_q$ (the blue solid line) and $b_q$ (the blue dotted line) 
when the NB terminated at iteration $q$ 
as $b_q < \hat{\zeta} = 441200000$ ({\it i.e.}, active node) --- Case (A) or 
at iteration $q$ as $\hat{\zeta} = 441200000 \leq a_q$ ({\it i.e.}, pruned node) --- Case (B).}
\vspace{-25mm}
\includegraphics[width=0.45\textwidth]{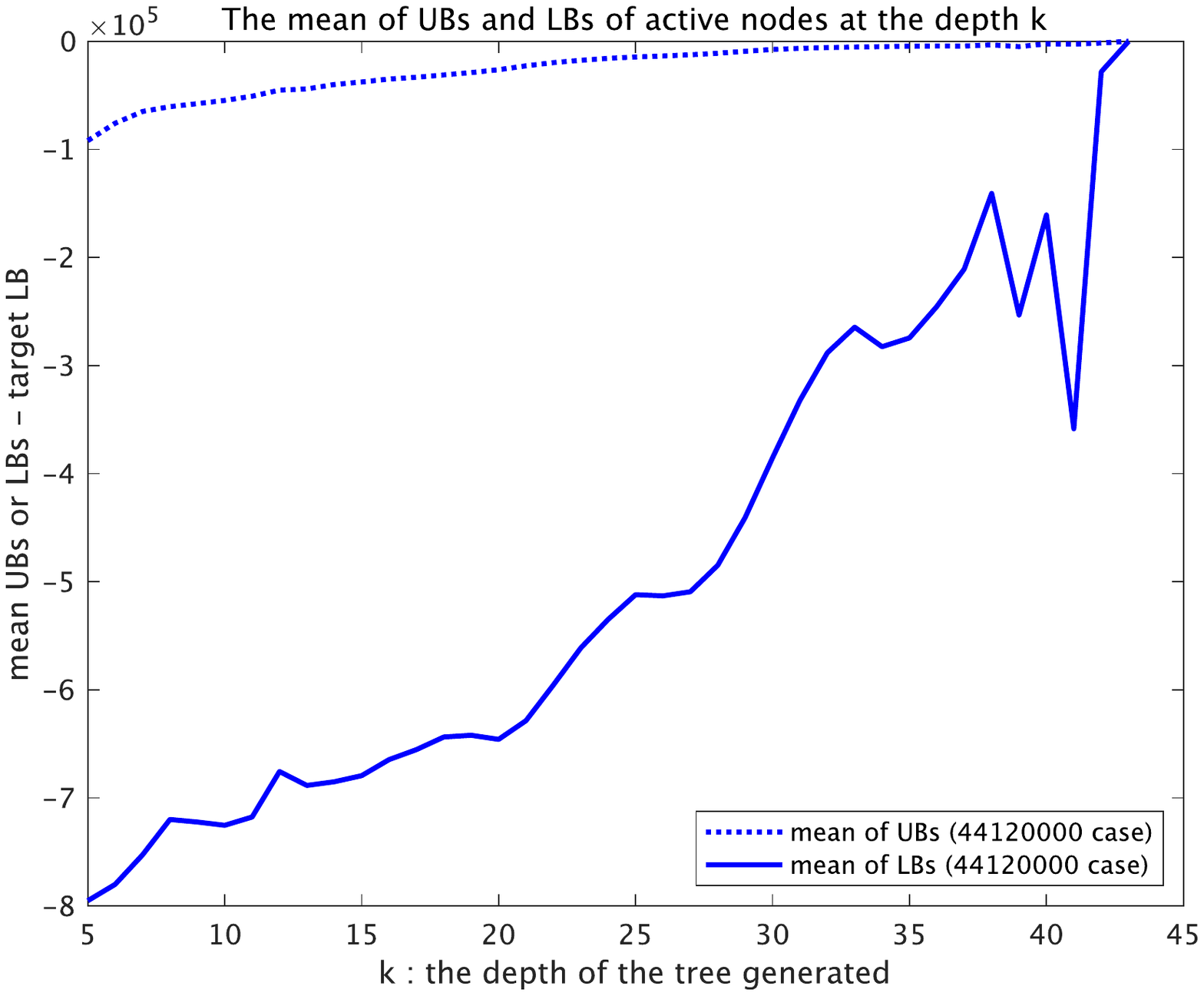} 
\includegraphics[width=0.45\textwidth]{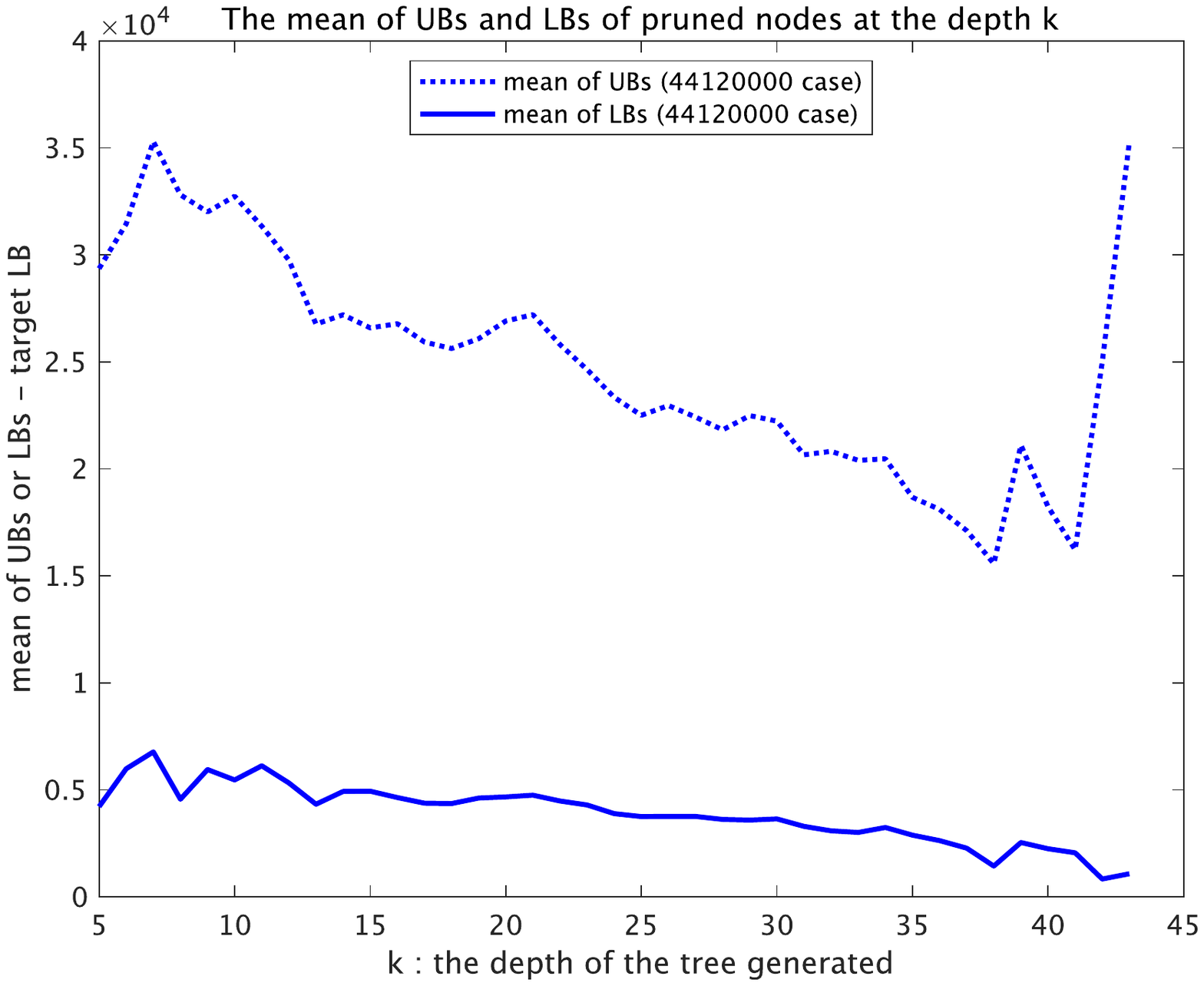} 
\end{center}
\end{figure}

\section{Concluding remarks}

\label{section:concludingRemarks}

While we 
have focused on the BQOP 
converted from the QAP tai256c in this paper, it is straightforward to adapt the discussion of the paper  
to  general BQOPs with a single cardinality 
constraint and  general QUBOs which satisfy the symmetry 
property~\eqref{eq:group0}.

We have shown that the largest unsolved QAP instance tai256c can be converted 
into a simple BQOP with a single cardinality constraint, and further into a 256-dimensional QUBO with 
a penalty parameter $\lambda$ whose optimal value converges to the optimal value of tai256c 
monotonically from below as $\lambda \rightarrow \infty$.  
Since the converted QUBO with dimension $256$  is much simpler than the original tai256c that involves $256 \times 256 = 65536$ binary variables, 
and its dimension $256$ is not so large, 
it might be natural to expect to solve it much easier than tai256c. We have found, however, that 
the QUBO is still quite difficult to solve. 
The difficulty essentially arises from 
the symmetry property~\eqref{eq:symmetry3} on the coefficient matrix $\B$ 
inherited from tai256c. We recall that 
the best known solution of tai256c is expanded to $1024$ different feasible solutions of the QUBO 
with the best known objective value. When  a branch and bound method is applied,  
all subproblems involving either of those feasible solutions need to be processed.


\end{document}